\documentclass[aps,pra,twocolumn,showpacs,superscriptaddress,nofootinbib]{revtex4}
\usepackage[dvips]{epsfig}
\usepackage{amsmath}
\usepackage{amssymb}
\usepackage{appendix}

 \graphicspath{{./}}

\newcommand{\real} {I\!\!R}

\newcommand{\gamovertwo} {{\frac{\gamma^2}{2}}}
\newcommand{\gammovertwo} \gamovertwo
\newcommand{\gammuovertwo} \gamovertwo
\newcommand{\gammuepstovertwo} \gamovertwo
\newcommand{\gammubarepstovertwo} \gamovertwo

\newcommand{\trace}{{\mbox{tr}}}

\newcommand{\midalign} {&&\hskip -1.5em}

\newcommand{\beasnum}{\begin{eqnarray}}
\newcommand{\eeasnum}{\end{eqnarray}}
\newcommand{\beas}{\begin{eqnarray*}}
\newcommand{\eeas}{\end{eqnarray*}}

\newcommand{\be}{\begin{equation}}
\newcommand{\ee}{\end{equation}}
\newcommand{\ba}{\begin{array}}
\newcommand{\ea}{\end{array}}

\newcommand{\figscale}{.35}

\def \numberofcontrolsinpracodfree{M}

\newtheorem{theorem}            {Theorem}[section]

\newtheorem{sideremark}         [theorem]{Remark}
\newtheorem{sideeg}           [theorem]{Example}
\newtheorem{sideconj}           [theorem]{Conjecture}

\newcommand{\qed} {\hskip 0.2em\lower 0.7ex\hbox{\vbox{\hrule
\hbox{\vrule height 1.2ex\hskip 0.4em\vrule height 1.2ex}
\hrule}}}

\newcommand{\SUn}{{SU(2^{n})}}
\newcommand{\lieg}{{\mathfrak{g}}}
\newcommand{\lieG}{\mathbf{G}}
\newcommand {\controlSet} {\mathcal{V}}

\newcommand{\ellfunction}[1] {\sqrt{{#1}^{T} R\,{#1}}\,}

\usepackage{color}

\begin{document}
\title{A reduced complexity numerical  method for  optimal gate synthesis}

\author{Srinivas Sridharan}
\affiliation{School of Engineering, Department of Information Engineering, Australian
    National University, Canberra, ACT 0200,
    Australia.}
\email[E-mail address to which correspondence should be sent: ]{srinivas.sridharan@anu.edu.au}
\author{Mile Gu}
\affiliation{Center for Quantum Technologies, National University of Singapore, Singapore}
\author{Matthew R. James}
\affiliation{School of Engineering, Department of Information Engineering, Australian
    National University, Canberra, ACT 0200,
    Australia.}
    \author{William M. McEneaney}
\affiliation{Department of Mechanical and Aerospace Engineering
		      University of California, San Diego, U.S.A}

\begin{abstract}
Although quantum computers have the potential to efficiently solve certain problems considered difficult by known classical approaches, the design of a quantum circuit  remains computationally difficult. It is known that the optimal gate design problem is equivalent to the solution of an associated optimal control problem, the solution to which is also  computationally
intensive. 
Hence, in this article, we introduce the application of a class  of numerical methods (termed the max-plus curse of dimensionality free techniques)   that determine the optimal control thereby synthesizing the desired unitary gate. The application of this technique to quantum systems has a growth in complexity that depends on the cardinality of the control set approximation rather than the much larger growth with respect to spatial dimensions  in  approaches based on gridding of the space, used in previous literature.  This technique is demonstrated by obtaining an approximate solution for the gate synthesis on $SU(4)$- a problem that is computationally intractable by grid based approaches. 

\end{abstract} 
\pacs{03.67.Lx, 02.70.-c}
\maketitle

\section{Introduction}

The advent of Shor\rq s algorithm \cite{shor1999polynomial} demonstrated the potential for processors based on quantum operations to perform certain computational tasks exponentially faster than those limited to to using classical operations. There has been much work devoted to solving the following problem of special interest: to determine the bounds on the number of one and two qubit gates required to perform a desired unitary operation--termed the gate complexity of the unitary.   Yet, the explicit design of quantum algorithms has remained a challenging task.

One approach  to this task of constructing an optimal circuit was highlighted in \cite{nielsen2006qcg} where it was shown to be  equivalent to finding a least path-length trajectory on a Riemannian manifold.  This insight opened up the study of quantum circuit complexity
to  the use of tools from optimal control theory.
 
In  \cite{sridharan2008gate}   the method of dynamic programming
  was introduced to solve the control problem associated with quantum circuit complexity. The numerical computations of solutions using this technique proceeded via a  widely used  grid (mesh) based  iteration approach \cite{kushner1992nms,Bardi,Crandall1984} that requires the generation of a mesh in the region of the state space over which the solution is sought. This approach however, leads to the following issue. A grid (assumed for simplicity to be regular and rectangular)  with $K$ points along each of the $N$ dimensions has $K^{N}$ grid points, over which the solution must be propagated during each iteration.  In addition, the dimension of an $n$ qubit quantum system  grows exponentially (as $4^n -1$) 
  thereby leading to a similar exponential growth in memory and time requirements. This large growth in the resources required,  arising from growth in the dimensions of the system, is termed  the curse of dimensionality (COD).   It renders the direct application of mesh based solution techniques unfeasible for systems larger than $SU(2)$ due to the large memory (in terabytes) and time (in centuries) required to solve  problems of these dimensions via such methods.

In  \cite{mceneaney2006mpm,mceneaney2008cdf} a COD-free technique was  introduced  for problems in Euclidean space. In this article we adapt these methods for quantum systems. Due to the structure of the control problem that we consider, we do  not completely eliminate the COD. However we have a much more manageable growth related to the number of elements in the discretized control set used. 
This is managed via a pruning approach described in Sec.~\ref{sec:pruning}. The computational time of the resulting algorithm grows much slower than  that in mesh based methods, thereby bringing us closer  to the numerical study of larger systems.  One particular application of interest for the numerical methods developed is the determination of  whether a given unitary $U$  in an $n$-qubit system can be  approximately synthesized in an efficient manner (with respect to the growth in $n$) in a given time $T$.

The paper is structured as follows. Sec.~\ref{sec:preliminaries} gives a brief introduction to the relevant concepts in quantum complexity and optimal control. We then introduce the reduced complexity algorithm in Sec.~\ref{sec:codfreemaxplus}, and in Sec.~\ref{sec:codfreemaxplusframeoptim} highlight the complexity growth in the application of this method and its management. The algorithm  is then applied  in Sec.~\ref{sec:exampleProblem}  to the two qubit optimal gate synthesis problem on $SU(4)$  - a problem in 15 dimensional space.  In Sec.~\ref{sec:pracodfree:conclusions} we conclude with comments on various aspects of the technique introduced in this article.

\section{Preliminary concepts} \label{sec:preliminaries}
In this section we recall the notion of gate complexity and introduce the  cost  function for an associated control problem as
in \cite{sridharan2008gate,Nielsen2006}.
\subsection{Gate complexity and control}\label{subsec:gateComplexity}
In quantum computing an algorithm operating on an $n$ qubit system can be represented as an element of the Lie group $\SUn$ (denoted in this article by $\lieG$) and is termed a unitary. Every such unitary can be constructed by a sequence of available elementary unitaries $U_{1}, U_{2} \ldots U_{n}$.  In practice,  we synthesize a unitary  $\hat{U}_0$ that approximates a desired computation $U_0$ with a required accuracy $\epsilon$ (i.e. $\|U_0-\hat{U}_0\| \leq \epsilon$, where
  $\|\cdot\|$ denotes the standard matrix norm).
This leads to the notion of approximate gate
complexity $G(U_0,\epsilon)$ which is 
  the minimal number of one and two qubit gates required to
  synthesize $U_0$ up to an accuracy of $\epsilon$ without ancilla qubits \cite{Nielsen2006}. 

  Related to the gate synthesis problem is an optimal control problem  (described below)  on $\lieG$, such that
the approximate  gate complexity  scales equivalently up to a polynomial in the optimal cost function for the control problem.
  This equivalence  motivates the solution of the associated control problem. 
  
  We now describe the control problem 
  and recall the solution process,  via the dynamic programming principle, used in  \cite{sridharan2008gate}.

\subsection{System Description} \label{sec:pracodfree:controlAndDP}
The system dynamics for the gate design problem is given by:
\begin{align}
\frac{dU}{dt} = -i\,\{\sum_{k=1}^ \numberofcontrolsinpracodfree v_k(t) H_k\} &U,\qquad
 U\,\in\,\lieG \label{eq:System}
\end{align}
with control $v$ (such that $v(t) \in \mathbb{R}^{\numberofcontrolsinpracodfree},\,\,\forall \, t\,\geq 0$) and an initial condition $U(0)=U_0$. 
For the class of problems  considered,  $v$ is taken to be an element of the set of piecewise continuous
functions having a norm bound $\|v(\cdot)\| =1$ (where $\| \cdot \|$ denotes the standard 2-norm on $\mathbb{R}^{\numberofcontrolsinpracodfree}$). We denote this class of controls by $\mathcal{V}$. The system equation contains a set
of right invariant vector fields $-i H_1,\,-i H_2\,\ldots\,, -i H_ \numberofcontrolsinpracodfree $, which
correspond to the set of available one and two qubit Hamiltonians. The
span of the set $\{ \, -i H_1,\,-i H_2\,\ldots\,, -i H_ \numberofcontrolsinpracodfree\}$  (and all brackets thereof) is
assumed to be the Lie algebra $\lieg$ of the group $\lieG$ \footnote{Note that we use the convention from mathematics  where elements of the Lie algebra are skew Hermitian. This is also consistent with the fact that the Hamiltonians are Hermitian.}.
Under these assumptions it follows from \cite[Prop. 3.15]{nijmeijer1990ndc} that the time to move the state, from the identity element to any other point on the group,  is bounded (and hence, the minimum time to move between any two points on $\lieG$ is finite). Given a control signal  $v$ and an initial unitary $U_0$ at time $r$ the solution to Eq~\eqref{eq:System} at time $t$ is denoted by
$U(t;v,r,U_0)$.

The control problem involves generating a desired state $U_{0}$ of the system in  Eq~\eqref{eq:System} starting from  the identity element. By time reversal of the dynamics Eq~\eqref{eq:System}  it can be seen that this is equivalent to the problem of reaching  the identity element starting from $U_0$.
The optimal cost function for this control problem is given by the geodesic distance
\begin{align}
  C_{0}(U_0) &= \mathop {\inf
  }\limits_{v\in{\controlSet_g}} \Big \{ { \int\limits_0^{t_{U_{0}}(v) } {\sqrt{{v(s)}^T R v(s)}}\,\, ds} \Big \} \label{eq:CostFunctionGeneral},\\
  \controlSet_g &:=\{v(\cdot)\,\big| \|v\|=1,\, \, \nonumber\\
&v_k\,:\,[0,\infty\,)\,\rightarrow\,\mathbb{R}\,\,\mathrm{\,is\,piecewise\,continuous} \nonumber\\ 
& \,\,\,\, \mathrm{\, for\, all \,k}\},\nonumber
\end{align}
where $t_{U_{0}}(v)$ is  the time to reach the identity starting from $U_{0}$ and is defined by
  \begin{align}
t_{U_0}(v) = \inf \{ t > 0 \ : \ U(0) = U_0, U(t)=I, \nonumber\\
\mathrm{dynamics\,in\,}  \mathrm{(\ref{eq:System})} \}.
\label{tUv-def}
\end{align}
This time is  taken to be $+\infty$ if the terminal
constraint $U(t)=I$ is not attained.

The  diagonal, symmetric and positive-definite weight matrix  $R$  in Eq~\eqref{eq:CostFunctionGeneral}   reflects the  relative difficulty of generating each element of the control vector. For instance,  on $SU(4)$, the  two body unitary  direction $\sigma_x \otimes \sigma_x$ may be weighted more than the single body unitary $\sigma_{x} \otimes I$ (as it is often harder to manipulate the former than the latter). The symbols $\sigma_{x}, \sigma_{y}, \sigma_{z}$ denote the standard Pauli matrices.

 One approach to optimal control problems such as in Eq~\eqref{eq:CostFunctionGeneral} is the \textit{dynamic programming} method developed in the 1950's by R.~Bellman (see \cite{bellman2003dp}). It relates the optimal cost evaluated at an initial time $r$ at a point $U_{0}$ to the optimal cost evaluated at a point $U(r+t;v,r,U_0)$ that is reached after applying the  control signal $v$ for a time duration of length $t$. This relation takes the form
%
\begin{align}
C_{r}(U_0) =&\mathop {\inf }\limits_{v\,\in\,\controlSet_{g}} \Bigg{\{}
\int\limits_r^{r+t }{{\sqrt{{v(s)}^T R v(s)}}\,\,ds} \quad + \nonumber\\ &  C_{r+t}(U(r+t ;v,r,U_0 ))
 \Bigg{\}}, \label{eq:DPEGeneral}
\end{align}
for all points $U_0\in \lieG$ and is termed the dynamic programming equation (DPE). For any function $\varphi$ this dynamic programming relation can be expressed as 
\begin{align}
C_{r}(U_{0}) &= S_{r,(r+t)}(C_{(r+t)})[U_{0}],\quad U_{0} \in \lieG \label{eq:pracodfree:dpo}
\end{align}
where 
\begin{align}
S_{r,(r+t)}(\varphi)[U_{0}]&:=  \mathop {\inf }\limits_{v\,\in\,\controlSet_{g}} \Bigg{\{} \int\limits_r^{r+t }{{\sqrt{{v(s)}^T R v(s)}}\,\,ds} \quad + \nonumber\\ &  \varphi(U(r+t ;v,r,U_0 ))
 \Bigg{\}}. \label{eq:Soperatordef}
\end{align}
Note that for ease of notation we will denote the operator  $S_{r,(r+t)}(\cdot)$ as $S_{t}(\cdot)$ when the value of $r$ is clear from the context.

In order to completely characterize the solution to the DPE \eqref{eq:DPEGeneral} we require boundary conditions given by 
\begin{align}
C(U) 
\left\{
\begin{array}{lr}
= 0,&U= I\\
>0,& U \in \lieG \setminus \{I\}.
\end{array}
\right.
 \label{eq:pracodfree:BC}
\end{align}
These conditions reflect the fact that   if we start at the identity $I$, then the cost to reach the identity is zero.  Furthermore, if $U_0 \neq I$, then a non-zero amount of time is needed to reach the identity as the  control values are bounded.
 
The solution (i.e. the optimal cost function) to the control problem can be obtained by solving a specific partial differential equation (a differential version of the DPE), termed  the Hamilton-Jacobi-Bellman (HJB) equation  \cite{Bardi}, given by
 \begin{align}
H(U,DC) &=  0, \quad U \in \lieG \setminus \{I\}\label{eq:pracodfree:HJB}
\end{align}
with
\begin{align}
 H(U,p) &:=\sup_{\|v\|=1 } \Bigg \{ {p \cdot \Big[ -i\,\big \{\sum_{k=1}^M v_k H_k \big \} \,U \Big]} \nonumber \\ &\qquad - \sqrt{{v}^T R v} \Bigg \}. \nonumber
\end{align}
 with  the boundary conditions  in Eq~\eqref{eq:pracodfree:BC}. 
The DPE/HJB  equation encodes considerable information about the problem, and can be solved for the optimal cost function. The DPE can then be used to construct/verify optimal control strategies via the \textit{verification theorem},   \cite[Sec.~1.5]{Bardi}).
 
 Various approaches exist to obtain the solution of the HJB equation \eqref{eq:pracodfree:HJB}. The most common are grid based methods \cite{Bardi, dupuis1999markov,fleming2006cmp, kushner1992nms} which require a mesh to be generated over the state  space. 
Due to the COD which leads to infeasible memory and time requirements these methods are unsuitable for larger dimension systems. Hence alternative approaches to this problem are required.

\section{The reduced complexity algorithm}\label{sec:codfreemaxplus}
Recently, a class of  algorithms  that are not subject to the curse of dimensionality were introduced in  \cite{mceneaney2008curse,mceneaney2008cdf,mceneaney2006mpm} to solve {first order} HJB equations in Euclidean space.
 This method, termed the  \textit{max-plus curse of dimensionality free} approach, to solve the dynamic programming equations involves a propagation of the solution of the HJB  equation forward in fixed time steps  without discretization in the spatial dimensions. The dramatic speed up of this approach stems from an invariant structure that  the cost function possesses, which is preserved under the above propagation. This invariant form helps reduce the amount of information that must be stored while solving the control problem. In this section we introduce this method and describe how the invariant structure arises from the cost function. We then apply this method to obtain a numerical procedure to determine the solution to the control problem.

The COD-free max-plus theory \cite{mceneaney2008cdf} to obtain and use an invariant form of the cost function  does not currently deal with cost functions containing terminal constraints  (as is the case in this problem where the trajectory  must reach the identity element  such as in Eq~\eqref{eq:CostFunctionGeneral}).  Hence  we  formulate a relaxed version of the problem, that can be solved via this  theory.

\subsection{Relaxation of the optimal cost function}
One possible relaxation of the cost function in Eq~\eqref{eq:CostFunctionGeneral} proceeds by
  introducing a fixed terminal time $T$ and a terminal penalty cost to yield the  expression
  \begin{align}
V^{\epsilon}_s(U_{0})&= \mathop {\inf }\limits_{v\in \controlSet^{e}_{g} } \Bigg\{ {\int\limits_s^{{T}} {\,{\sqrt{{v(t)}^T\, R \, v(t)}}\,\,dt} }  \,+ \nonumber \\  &{\,\,\frac{1}{\epsilon} \phi\big(U({T};v,s,U_{0})\big)} \Bigg\}, \label{def:costEpsilon}
\end{align}
where  $\phi(\cdot)$ is a real valued non-negative function (that is zero only at the identity element). This function  penalizes terminal states away from the identity.  

The extended control set $\controlSet^{e}_{g}$ above, is defined as 
\begin{align}
\controlSet^{e}_{g} : = \controlSet_{g} \,\, \bigcup \,\,\{v \equiv 0 \},
\end{align}
where $\{v \equiv 0 \}$ denotes the control signal that is  identically zero.
We note that due to the fixed time horizon in Eq~\eqref{def:costEpsilon} this extension to the control set ensures that once the target set is reached, the cost function does not increase
further.

 In this article we take the terminal cost to be of the form 
\begin{align}
\phi(U)&=  \trace[(I-U)\times (I-U)^{\dag}]  \nonumber \\&=  \trace[2I-U-U^\dag], \label{eq:pracodfreetermcostdef}
\\&= 2 \,\trace[I] - 2 \mathop{Re}(\trace(U)).
\end{align}
where \lq $\trace$\rq\, denotes the trace  operation and $\mathop{Re} (\cdot)$ denotes the projection onto the real axis.
The above relaxation is valid since, as the penalty increases, the aproximation $V_{0}^{\epsilon}(U)$ converges to $V(U)$ for all points $U \in \lieG$ (for a sufficiently large time horizon $T$).

  We recall that  due to the assumptions outlined in Sec.~\ref{sec:pracodfree:controlAndDP},  the time to move between any 2 points in the group $\lieG$ is bounded.  The minimum time to move from $U_{0}$ to $I$ is denoted by
  \begin{align}
t_{U_0} = \inf_{v \in \mathcal{V}^{e}_{g}}\big[t_{U_0}(v) \big] .
\label{tU-def}
\end{align}
We define
\begin{align}
{T_{0}}:= \sup_{U_{0} \in \lieG} [t_{U_0}] \quad (< \infty)
\end{align}
to be the maximum value of this time over all points in the group. Hence choosing a  ${T}$ to be $T_{0}$ in  Eq~\eqref{def:costEpsilon} ensures that it is applicable for any  initial point $U_{0} \in \lieG$.

We now propagate the relaxation $V^{\epsilon}_s$, that  satisfies the DPE in Eq~\eqref{eq:pracodfree:dpo}, in time steps of $\tau$.
%
Let $T~=~{N} \tau$, for some ${N} \in \mathbb{N}$.   Using the notation 
$$
\tilde{V}^{\epsilon}_k(U) = V^{\epsilon}_{T-k \tau}(U)
$$
and  Eq~\eqref{eq:pracodfree:dpo}  we  rewrite Eq~\eqref{def:costEpsilon} as  the propagation of  the function $\tilde{V}^{\epsilon}$ forward in time through  a time interval $\tau$  due to  the action of the operator $\bar{S}_{\tau}[\cdot]$:
\begin{align}
\tilde{V}^{\epsilon}_{k+1}(U) = \bar{S}_{\tau}[ \tilde{V}^{\epsilon}_k(\cdot)] (U),\label{eq:pracodfree:soperatorpropagation}
\end{align}
where $\bar{S}_{\tau}$ is the form of the operator $S_{\tau}$ which uses piecewise constant controls over each time step of duration $\tau$.

Hence we  repeatedly apply $\bar{S}_\tau[\cdot]$ to move towards the desired value function $V^{\epsilon}_0(\cdot)$ ( $\cong  \tilde{V}^{\epsilon}_{N}(\cdot)$).

From Eq~\eqref{def:costEpsilon}, \eqref{eq:pracodfree:soperatorpropagation}  it follows that for all $\tilde U \in \lieG$,
\begin{align}
\tilde{V}^{\epsilon}_{k+1} (\tilde U)=&\mathop {\inf }\limits_{v\,\in\,\controlSet^{e}_{\tau}} \Big \{
(\ellfunction{v})\,\times \tau +  \nonumber \\ &\tilde{V}^{\epsilon}_{k} (U(\tau ;v,k \tau, \tilde U))\Big \}, \label{eq:pracodfree:simplerproptau}
\end{align}
which indicates the action of the  DPE operator $\bar{S}_{k \tau,(k+1)\tau}[\cdot]$. We now describe how the properties of this operator  lead to  the  invariant structure of the optimal cost function Eq~\eqref{def:costEpsilon}, and the implications resulting therefrom.

\subsection{Invariant structure of the optimal cost function} \label{sec:codfreemaxplus:invariant}

The terminal cost  in Eq~\eqref{eq:pracodfreetermcostdef} can be written as
\begin{align}
\phi(U) = c_{0} + P_{0}(U),\, U \in \lieG, \label{eq:pracodfree:originalP}
\end{align}
where 
\begin{align}
c_{0} = 2 n, \qquad P_{0}(U) = - \trace [U] - \trace[{U}^{\dag}].
\end{align}
Hence $c_{0}, P_{0}$ encode the initial costs corresponding to a control $I$ (i.e. no control action).   Assume that for a given \lq $k$\rq\, ($k \in [1, 2, \ldots N-1]$), the cost function $\tilde V^{\epsilon}_{k}$ can be written as 
\begin{align}
\tilde V^{\epsilon}_{k} (U) = c_{k} + P_{k} (U).\label{eq:pracodfree:vinitialpreservation}
\end{align}
From Eqns~\eqref{eq:pracodfree:soperatorpropagation} and \eqref{eq:pracodfree:simplerproptau}, after one time step $\tau$ the value function becomes  
\begin{align}
\tilde V^{\epsilon}_{k+1}(U_{0}) &= \bar{S}_{\tau}[\tilde V^{\epsilon}_{k}](U_{0}) \\
&= \min_{v \in \controlSet^{e}_{\tau}} \{ (\ellfunction{v})\, \tau + c_{k} + \nonumber \\ &\quad P_{k}(U((k+1)\tau;v,k \tau, U_{0})) \},\\
&=  \min_{v \in \controlSet^{e}_{\tau}} \{(\ellfunction{v})\, \tau + c_{k} + \nonumber \\ & \quad P_{k}(\Psi[(k+1) \tau, k \tau, v] \cdot U_{0}) \},\label{eq:pracodfree:plam}\\
&:=\min_{v \in \controlSet^{e}_{\tau}} p^{v}_{k}(U_{0}), \label{eq:pracodfree:useofp}
\end{align}
where $\forall U_{0} \in \lieG$
\begin{align}
p^{v}_{k}(U_{0}) &:= \{(\ellfunction{v})\, \tau + c_{k} +\nonumber\\& P_{k}(\Psi[(k+1) \tau, k \tau, v] \cdot U_{0}) ,  \label{eq:pracodfree:defofp}\\
\Psi[t, s, v] \cdot U_{0}&:= U(t;v,s,U_{0}).
\end{align}
In the equations above, $\Psi(\cdot)$ denotes the propagator for the system dynamics in Eq~\eqref{eq:System}  under the action of a  control signal $v$.

Hence 
\begin{align}
\tilde V^{\epsilon}_{k+1}(U_{0}) &= c_{k+1} + P_{k+1}(U_{0}), \label{eq:pracodfreeP1andP0}
\end{align}
where
\begin{align}
 c_{k+1}&:= c_{k} + (\ellfunction{\bar v}) \tau,\\
P_{k+1}(U_{0}) &:=  P_{k}(\Psi[(k+1) \tau, k \tau,\bar v] \cdot U_{0}),\\
\text{where} \quad \bar v &=\mathop{\arg\min}_{v \in \controlSet^{e}_{\tau}} \big [ p^{v}_{k}(U_{0}) \big] .
\end{align}

By the principle of induction, from  Eqns.~\eqref{eq:pracodfree:originalP}, \eqref{eq:pracodfree:vinitialpreservation}, \eqref{eq:pracodfreeP1andP0}  it may be seen that  $\bar{S}_{\tau}[\cdot]$ preserves the structure of the cost function. This invariance of the structure is a key aspect of the class of techniques introduced, as it helps obtain the optimal  cost function  at desired points without having to discretize along  the spatial dimensions. This optimal cost function for any point $U_{0}$ is 
\begin{align}
\tilde{V}^{\epsilon}_{N}(U_{0}) = \min \limits_{k \in \{1, 2 \ldots N\}} \min_{v \in\, {\prod \limits_{k}\controlSet^{e}_{\tau}}} p^{v}_{k} (U_{0}),\label{eq:pracodfree:gettingcostfrmp}
\end{align}
where $\prod \limits_{k}$ denotes the $k$ fold product of the control set $\controlSet^{e}_{\tau}$.
Hence once a computationally efficient parameterization of the control signals in terms of the set of $p$ values is obtained as described above, Eq~\eqref{eq:pracodfree:gettingcostfrmp} easily yields the cost function. We note that the computation of $P_{k}$ and $ \Psi[t,s,v] \cdot U_{0}$  can be performed  efficiently as they can be reduced to matrix multiplications and trace operations on matrices.  The generation of a set of parameters, by using the invariant structure of the cost, and its application to determine the optimal cost function is the essence of the max-plus COD free technique. 

From Eqns~\eqref{eq:pracodfree:useofp},\eqref{eq:pracodfree:defofp}  it is clear that during each time step there is an increase in the size of the number of candidate controls to be considered during the minimization. Specifically, the number  of elements $p^{v}_{k}[\cdot]$ that result from each  $P_{k}$ is the cardinality of the control set $\controlSet^{e}_{\tau}$. Thus, as in the COD free method in Euclidean space,  due to  the avoidance of spatial discretization the problem is free of the growth in dimensionality arising from spatial terms; however due to the structure of the quantum control problem \textit{there is now a geometric growth in complexity}. The details of this growth and methods to reduce its impact are now described.

%
\section{Control space growth and pruning} \label{sec:codfreemaxplusframeoptim} \label{sec:pruning}

For the purpose of implementation let the control space be discretized as follows: the control signal is held  constant over  each particular time period of duration $\tau$. 
Furthermore atmost one component of the control vector  is set to a value of $1$ over any time period, while the others are kept at $0$.  This class of control signals is denoted by $\tilde \controlSet^{e}_{\tau}$.   
Note from  Eq~\eqref{eq:pracodfree:plam} that  after each time step there is a factor of $[\#(\tilde \controlSet^{e}_{\tau})]$  growth in the number of control sequences in the set to be considered, where $\#(A)$ indicates the cardinality of a set $A$. Hence after   $N$ time steps the number of possible control sequences is $[\#(\tilde \controlSet^{e}_{\tau})]^{N}$. 

To manage this growth we introduce a selective removal (termed pruning) of some of these control sequences
 To describe this pruning procedure we first introduce the required  notation. The set of control sequences of length $k$  (i.e.  $k$ time step sequence) is denoted  by $\Lambda_{k}$. 
The set of all such control sequences of all possible lengths $\{1,2...\ldots N\}$ is 
$$
\Lambda:= \bigcup_{k=1}^{N} \Lambda_{k}.
$$

As indicated in the previous section there is a cost function $p_{\lambda}(\cdot)$, associated with each control $\lambda \in \Lambda$ such that
the discretized cost function $\tilde{V}^{\epsilon}_{N}(U)~\cong~ V^{\cdot}_0(\cdot)$  can be determined for any point $U \in \lieG$ by
\begin{align}
\tilde{V}^{\epsilon}_{N}(U) & =  \min_{\lambda\in\Lambda} p_{\lambda}(U), \label{eq:pracodfree:totalminimization}
\end{align}

To decide upon pruning some of the control sequences $p_{\lambda}(\cdot)$ we first determine its contribution to the minimization of the cost function. A control $\bar \lambda$ (and the corresponding function $p_{\bar \lambda}$) contributes to the minimization of the cost function Eq~\eqref{eq:pracodfree:totalminimization} iff
\begin{align}
\exists U \in \lieG\,\, \text{such\,that}\,\, p_{\bar \lambda}(U) < p_{\lambda}(U),\quad \forall \lambda \in \bar \Lambda,
\end{align}
where $\bar\Lambda$ denotes the  set of control sequences ${\Lambda\setminus \{\bar\lambda\}}$ which are different from the sequence  $\bar{\lambda}$.
This idea can be used to measure the contribution of any control sequence towards the minimization of the cost in Eq~\eqref{eq:pracodfree:totalminimization}.  From \cite{mceneaney2008curse} one such function that quantifies  this contribution  is 
\beasnum
W(\bar\lambda)\doteq\max\bigl\{h(\zeta,U)\,\big\vert\,\midalign\bar p_\lambda(U)-\zeta\ge 0,\,\forall \lambda\in\bar\Lambda,\,\nonumber\\
U \in \lieG, \,\zeta\in\real\bigr\},
\label{eq:W1}\eeasnum
where 
\begin{align}
\bar p_\lambda(U) &\doteq p_{\lambda}(U)-p_{\bar\lambda}(U)\quad \forall \lambda \in \bar\Lambda, \nonumber\\
h(\zeta,U)&\doteq \zeta, \quad \zeta\in\real,
\end{align}

If $W(\bar\lambda)\le 0$, then $p_{\bar\lambda}$ never achieves the minimum for any  point $U$ in 
$\lieG$, and consequently, the control sequence ${\bar\lambda}$ can be pruned without any effect on  the cost function.
%
For those $\bar\lambda$ such that $W(\bar\lambda)>0$, pruning would remove control sequences that do contribute to the minimization. However this is unavoidable in order to manage the growth in computational resources required. Therefore  to reduce the errors in the optimal cost function arising from this pruning, we  selectively eliminate control sequences with 
relatively small values of $W$. This minimizes the impact of pruning on the optimality of the resulting control strategy. 

There are several numerical methods \cite{ben2001lectures,YALMIP} that can efficiently solve  pruning problems of the form in Eq~\eqref{eq:W1}.   The more involved mathematical details of the procedure will be addressed in a subsequent article.
By adjusting the upper limit on the number of control sequences stored in each set $\Lambda_k$, we may arrive at an acceptable tradeoff between speed and accuracy. 
This approach has enabled a dramatic improvement in the time required to solve problems with 2-qubits (Sec.~\ref{sec:exampleProblem}), while using  standard computing resources

\subsection{Description of computational complexity of the algorithm} \label{sec:pracodfree:compcomp}
We now compare the complexity of the algorithm outlined in this article with that  of mesh based solution methods such as in \cite{sridharan2008gate}.  In the reduced complexity method without pruning, the computational complexity grows as
\begin{align}
\mathcal{M}^{K} \label{eq:pra:codfree:compcomp}
\end{align}
 where $\mathcal{M}$ is the number of elements in the control set and $K$ is the number of time steps in the simulation. For the assumptions on controllability to hold, at-most $\mathcal{M}$ ($\in O(n^{2})$)  directions of control  (i.e. control Hamiltonians) are required. 
  Hence, from Eq~\eqref{eq:pra:codfree:compcomp},   the complexity of the algorithm (without pruning) for a simulation of $K$ time steps  is 
\begin{align}
[O(n^{2})]^{K}  = \mathop{poly} (n),\label{eq:pracodfreecomplexitywithpruning}
\end{align}  
where $\mathop{poly} (n)$ denotes a polynomial in $n$.
With pruning, the complexity growth depends on the storage limits chosen in  the pruning process. Hence there exists a tradeoff, influenced by these storage limits, between the accuracy of the solution and the growth in  complexity of the  procedure required to obtain it.

In contrast  the computational cost in Eq~\eqref{eq:pracodfreecomplexitywithpruning}, for mesh based methods  with $\Gamma$ mesh points along each dimension, is 
$$
\Gamma^{[4^{n}-1]} \times  O(\Gamma).
$$
The two terms in this expression arise from the number of spatial dimensions and number of iterations required respectively.

 An important application  of the approach described herein is that, {given a fixed time horizon $T$  it is possible to efficiently (with respect to the scaling of $n$) check if a desired gate $U$ can be synthesized within this time}. The complexity of the algorithm to perform this check  this would be $[O(n^{2})]^{K}$ (\textit{without pruning}) where $K$ is the number of iterations in the algorithm  (which is fixed for a given $T$).

Thus it may be observed that the COD-free approximation  technique offers a potentially large order of magnitude  reduction in computational complexity.

\section{Example on $SU(4)$}\label{sec:exampleProblem}
We now proceed to apply the theory introduced, to an example on  $SU(4)$. 
The dynamics for this system is given by Eq~\eqref{eq:System}  with $\numberofcontrolsinpracodfree =5$ and a  control set  generated from Hamiltonians of the form
 $H_k \,\in\,\{I\otimes \sigma_x, \, I \otimes \sigma_z,\,\sigma_x\otimes I,\,\sigma_z\otimes I, \sigma_x \otimes \sigma_z\,\}$ i.e., a set of four $1$-body terms and one $2$-body term. The associated control directions are sufficient to generate the entire Lie algebra $\mathfrak{su}(4)$, thereby ensuring controllability.

To help highlight  the performance improvements of the methods introduced herein, we note that  a grid based solution approach with a conservative mesh of $50$ points in each of the $15$ dimensions  of the space $SU(4)$, with a total of $20$ iterations over the space and  an estimated time of $0.0001$ seconds to propagate the cost function via value iteration \cite{sridharan2008gate} at each point in the mesh would require $1.69\times 10^{19}$ hours and a few terabytes of memory.
Using the reduced complexity theory described in the previous section, this problem was  solved in $15$ hours to yield a solution  for the final time horizon problem, with a horizon $T$ of  $4$ seconds and a discretization step size $\tau$ of $0.2$ seconds (i.e., 20 propagation steps). The simulation was carried out on a  standard desktop computer without any exhaustive efforts to optimize the code. Hence there is a strong potential for further improvements to this procedure. 
\subsection{Simulation results}\label{sec:pracodfree:simresults}
The simulation results obtained from the reduced complexity technique provide the optimal cost function Eq~\eqref{eq:pracodfree:totalminimization} for the gate synthesis problem on the two qubit system. 
In order to visualize the cost function on the group, we  require a mapping between points on the group and points in Euclidean space (as the latter can be easily plotted via conventional graphs). For this purpose we make use of the exponential map from the literature on differential geometry  \cite{Sepanski2007,hall2003lgl}. This is an onto map that takes points in the 
Lie algebra (the tangent space at the identity element) to points in the group. As the Lie algebra $\lieg$ is isomorphic  to the Euclidean space, we can thus obtain a function in the Euclidean co-ordinates at points  of interest.  The exponential map acts on the algebra of  any matrix Lie group   as follows
\begin{align}
\exp(X) := \sum_{j=0}^{\infty} \frac{X^{j}}{j!}, \quad \forall X \in \lieg.
\end{align}
For the connected group $SU(2^{n})$  the exponential map of the algebra ($\exp(\lieg)$) generates all of the group ($\lieG$)\cite[Thm 4.6]{Sepanski2007}, thereby ensuring that a valid visualization can be generated for all points in the group.

\begin{figure}[thp]
\centering
\includegraphics[scale= \figscale]{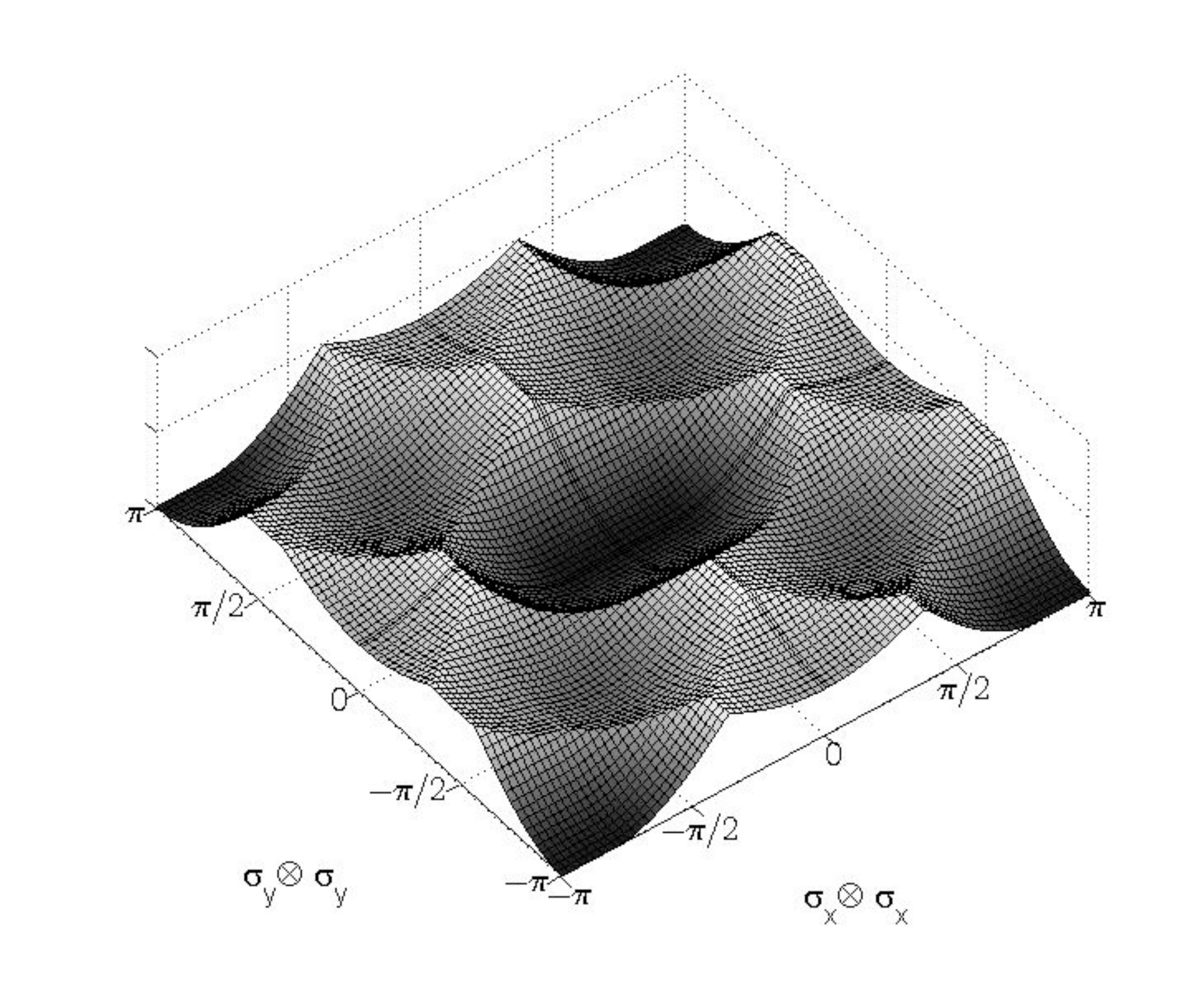}
\caption{Plot of the cost function in the plane $\sigma_x~\otimes~\sigma_x$ vs $\sigma_y~\otimes~\sigma_y$. Here the cost  of applying any available 1-body Hamiltonian is taken to be only slightly less than   the cost of generating the 2 body Hamiltonian $\sigma_x~\otimes~\sigma_x$ (i.e., $\mathbf{r}$ = 1/1.3). The lighter areas indicate unitary operations that are harder to synthesize.}
 \label{fig:pracodfree:figixxiyycost1p3}
\end{figure}

Thus, in order to visualize the cost function on $\lieG$, we evaluate the optimal cost function at a desired set of points on this group that correspond to a 2-dimensional slice of interest in $\lieg$. %
The plots obtained  indicate the approximate optimal cost function at points chosen  in the plane of interest (on the algebra)  e.g., the $\sigma_x \otimes \sigma_x$ vs $\sigma_y \otimes \sigma_y$  plane (used in this article). 
The value of the cost function is mapped to the shading used, in order to illuminate the behavior of the function. For instance in Fig~\ref{fig:pracodfree:figixxiyycost1p3}, regions of darker shading indicate unitaries which are easier (lower cost) to generate while the lighter areas 
show gates which are costly to synthesize.  

\begin{figure}[thp]
\centering
\includegraphics[scale= .3]{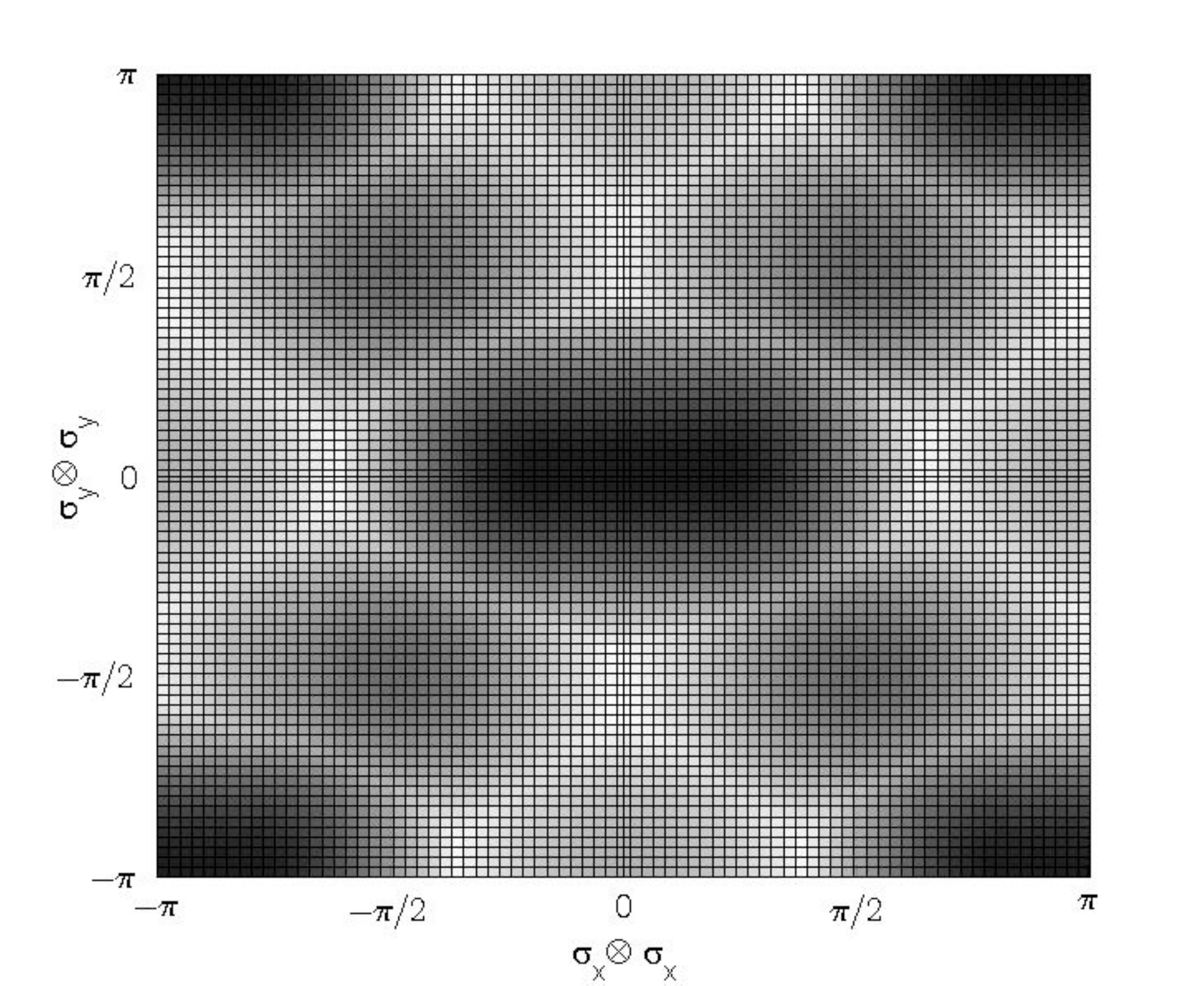}
\caption{Top view of the plot in figure \ref{fig:pracodfree:figixxiyycost1p3}.}
 \label{fig:pracodfree:figixxiyycost1p3flat}
\end{figure}

In order to understand the effects, on the cost function, of the difficulty  in synthesizing the $2$-body unitary compared to the $1$-body terms (which we denote by the ratio $\mathbf{r}$ of the cost of one body interactions  to that of the available two body interaction $\sigma_x \otimes \sigma_x$ ), simulations were performed that varied $\mathbf{r}$ (using the $R$ matrix in Eq~\eqref{eq:CostFunctionGeneral}) starting from  slightly less than one  and proceeding upto a value of  one-third.

\begin{figure}[h!]
\centering
\includegraphics[scale= \figscale]{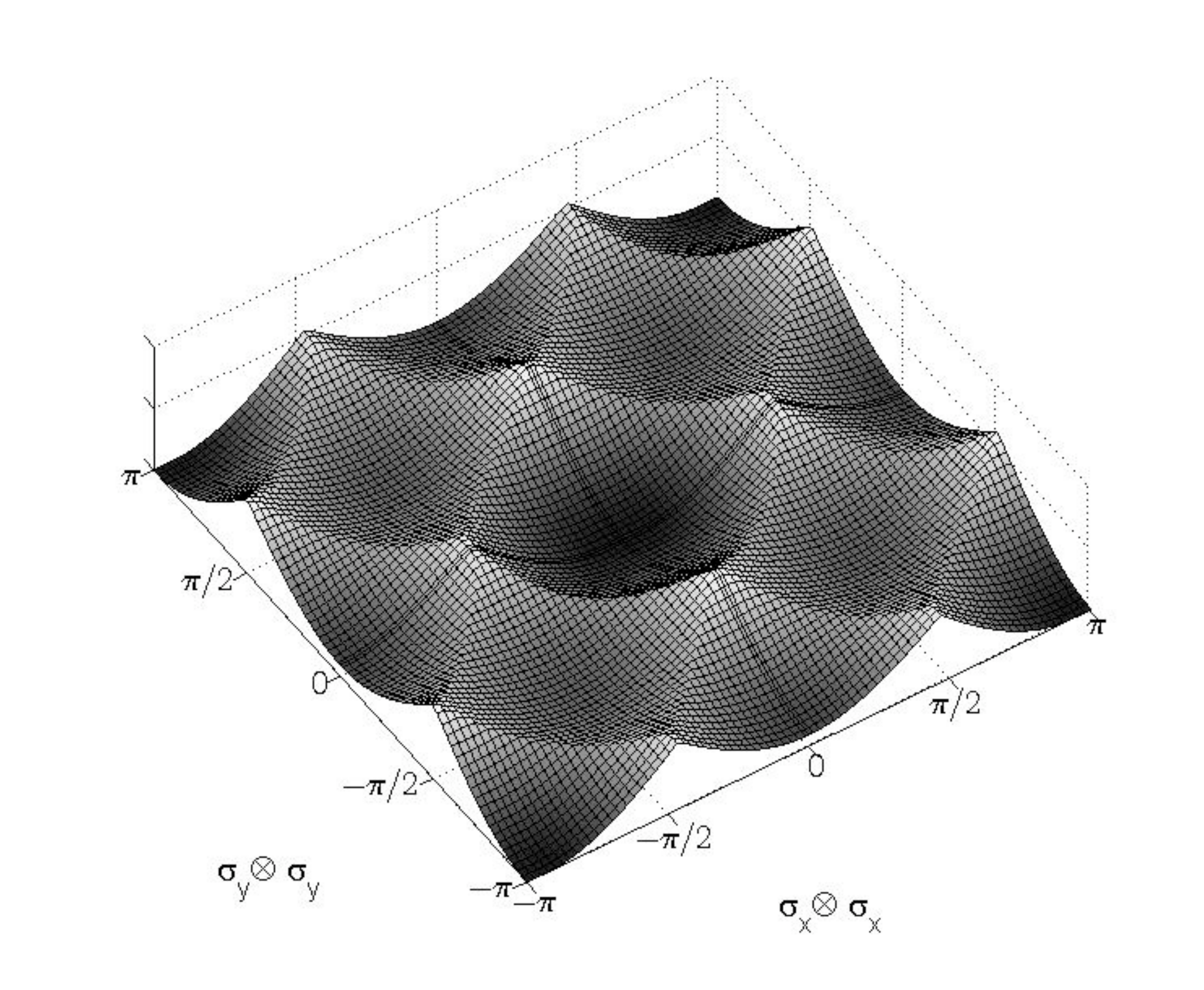}
\caption{Plot of the cost function in the plane  $\sigma_x \otimes \sigma_x$ vs $\sigma_y \otimes \sigma_y$ with $\mathbf{r}$   set to $1/3$.}
\label{fig:pracodfree:figixxiyycost3}
\end{figure}

\begin{figure}[h!]
\centering
\includegraphics[scale= .3]{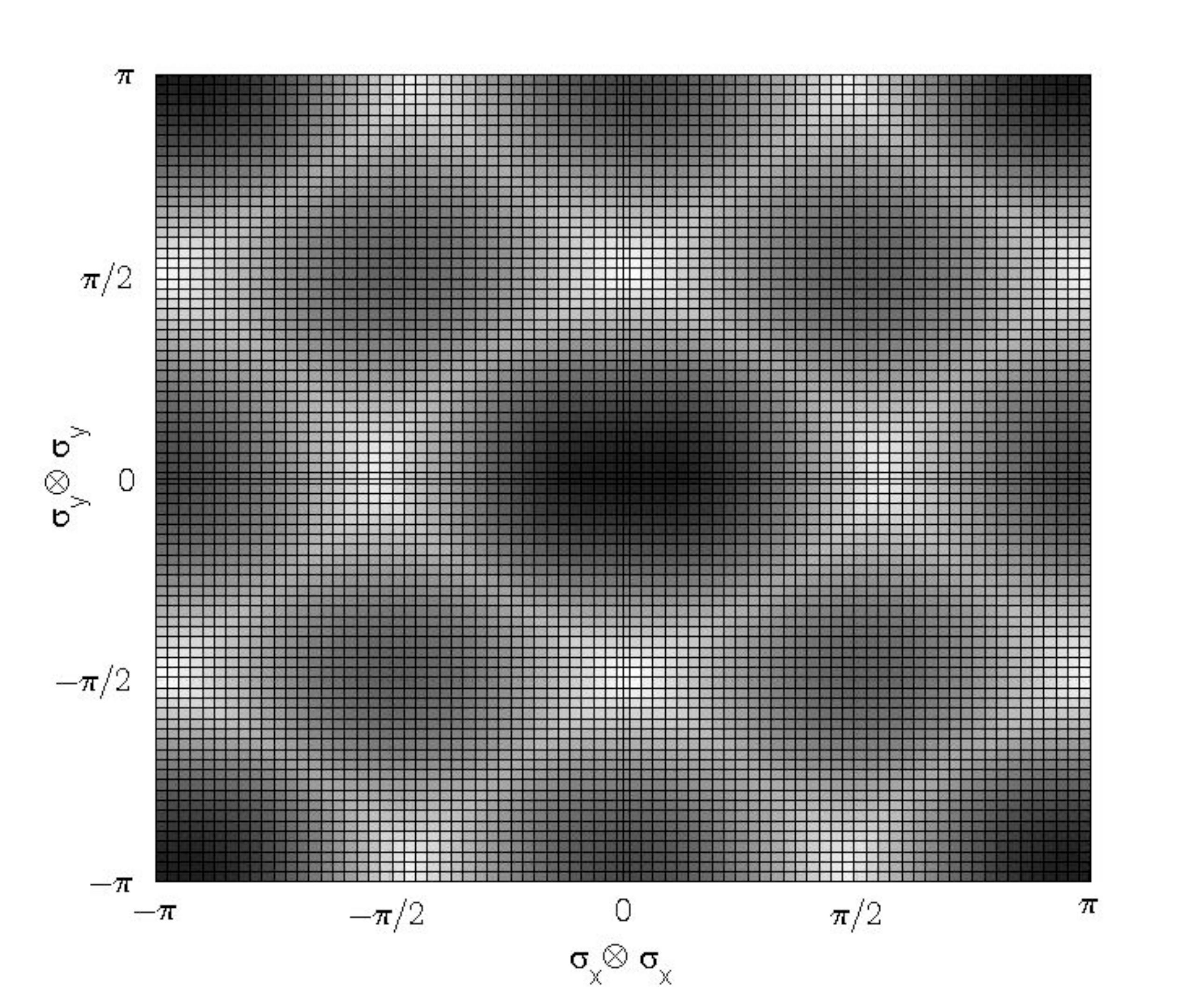}
\caption{Top view of the plot in figure \ref{fig:pracodfree:figixxiyycost3}.}
\label{fig:pracodfree:figixxiyycost3flatinterp}
\end{figure}

The generation of unitaries in the $\sigma_y\otimes \sigma_y$ direction (which is not directly accessible) requires the  alternating  application  (i.e., bracketing operation) of 
$\sigma_x\otimes \sigma_x$ and one of the  available single body Hamiltonians. Hence if $\mathbf{r}=1/1.3$,  then the cost of moving along the $\sigma_y \otimes \sigma_y$  direction is substantially larger than that along the one and two body control directions, leading to an elliptical shape of the level set (as shown in Fig.~\ref{fig:pracodfree:figixxiyycost1p3flat}), where the cost increases faster in the $\sigma_y \otimes \sigma_y$ direction compared to the $\sigma_x \otimes \sigma_x$ direction.  However as the cost  imposed on moving along the $\sigma_x \otimes \sigma_x$ direction is increased, the relative cost of moving along $\sigma_y \otimes \sigma_y$  compared to $\sigma_x \otimes \sigma_x$  decreases, thereby leading to level sets that become more circular (Fig.~\ref{fig:pracodfree:figixxiyycost3flatinterp}). This agrees with analytical results such as in \cite{phyrevGuexplicitexample}. Note that  there are two axis of symmetry in these plots, namely both $\sigma_x\otimes \sigma_x$  and $\sigma_y\otimes \sigma_y$. This arises due to the symmetry in the cost function equations \eqref{eq:CostFunctionGeneral} and  \eqref{def:costEpsilon} with respect to the application of  controls along either the positive or negative direction of the available control Hamiltonians in the system equation \eqref{eq:System}.\\

\section{Conclusions}\label{sec:pracodfree:conclusions}
In this article we have  demonstrated a reduced complexity method may be used to obtain substantial improvements in computational speed in solving optimal control problems arising in  closed quantum systems. This technique was used to obtain a numerical solution for an example gate complexity problem in a $2$ qubit system.  
Instead of the curse of dimensionality in the spatial dimension, we now have a much more manageable growth in dimensionality  that depends on the number of elements in the discretization of the control set.

The approach outlined in this article can deal with a very general class of problems and can, in principle, be used for systems with any number of spins.  Furthermore, the methods described  can be extended to other systems of interest. For instance, a control problem on a system with drift can be solved under the current system framework (Eq~\eqref{eq:System}) by taking the cost of moving along the negative direction of the drift term to be much larger than that in the positive direction. 

At present, the techniques outlined yield preliminary solutions whose error bounds, rate of convergence and other properties must be determined via further research. It is hoped that the  method introduced (and refinements thereof) would enable the accurate solution of problems on spin systems of  larger dimensions than has been possible until now.

\begin{acknowledgments}
S.~Sridharan and M.R.~James wish to acknowledge the support for this work by the
    Australian Research Council. M.~Gu acknowledges the support from the National Research Foundation and Ministry of Education of Singapore.  W.~McEneaney acknowledges  support from  AFOSR and NSF.  The authors would like to thank the reviewer for helpful comments.
    \end{acknowledgments}
\bibliographystyle{apsrev}

\end{document}